\documentclass[12pt]{amsart}
\usepackage{diagram}
\usepackage{latexsym}
\input amssym.def
\input amssym.tex
\numberwithin{equation}{section}
\makeatletter
\renewcommand{\subsection}{\@startsection
{subsection}{2}{0mm}{\baselineskip}{-0.25cm}
{\normalfont\normalsize\rmfamily}}
\makeatother
\newtheorem{claim}{Claim}
\newtheorem{theorem}{Theorem}[section]

\newtheorem{corollary}[theorem]{Corollary}
\newtheorem{prop1}{Proposition}[subsection]
\newtheorem{lemma1}[prop1]{Lemma}
\newtheorem{corollary1}[prop1]{Corollary}
\newtheorem{thm1}[prop1]{Theorem}
\theoremstyle{remark}
\newtheorem{remark}[prop1]{Remark}
\newenvironment{claim*}{\begin{trivlist}\item[\hskip%
\labelsep{\bf{Claim.}}]\it }%
{\end{trivlist}}
\newtheorem{remarke}[theorem]{Remark}
\def\d{{\mathcal D}}
\def\j{{\mathcal J}}
\def\T{T_{\ell}({\mathcal J})}
\def\pij{{\rm Fr}_{\mathcal J}}
\def\pix{{\rm Fr}_{X}}
\def\supp{{\rm Supp}}
\def\div{{\rm div}_\infty}
\def\h{{\mathcal H}_{m,q}}
\def\kx{{k}(X)}
\begin{document}
\author[R.~Fuhrmann]{Rainer Fuhrmann}\thanks{The paper was written while
the first author was visiting the Instituto de Matem\'atica Pura e
Aplicada, Rio de Janeiro (Oct. 1995 - Jan. 1996). This visit was
supported by Cnpq}
\author[F.~Torres]{Fernando Torres}\thanks{The second author is supported
by a grant from the International Atomic Energy Agency and UNESCO}
\title[Curves with many rational points]{On curves over finite fields
with\\ many rational points}
\address{Fachbereich 6, Mathematik und Informatik, Universit\"at Essen,
D-45117 Essen - Germany}
\email{100713.3165@compuserve.com}
\address{ICTP, Mathematics Section, P.O. Box 586, 34100, Trieste - Italy}
\email{feto@ictp.trieste.it}
\begin{abstract}
We study arithmetical and geometrical properties of {\it maximal curves},
that is, curves defined over the finite field $\mathbb F_{q^2}$ whose
number of $\mathbb F_{q^2}$-rational points reachs the Hasse-Weil upper
bound. Under a hypothesis on non-gaps at rational points we prove that
maximal curves are $\mathbb F_{q^2}$-isomorphic to $y^q+y=x^m$ for some
$m\in \mathbb Z^+$.
\end{abstract}
\maketitle
\setcounter{section}{-1}
\section{Introduction}
Goppa in \cite{Go} showed how to construct linear codes from curves
defined over finite fields. One of the main features of these codes is
the fact that one can state a lower bound for the minimum distance of the
codes. In fact, let $C_X(D,G)$ be a Goppa code defined over a curve $X$
over the finite field $\mathbb F_q$ with $q$ elements, where
$D=P_1+\ldots+P_n$, $P_i\in X(\mathbb F_q)$ for each $i$ and $G$ is a
$\mathbb F_q$-rational divisor on $X$. Then it is known
that the minimum distance $d$ of $C_X(D,G)$ satisfies
$$
d\ge n - {\rm deg}(G).
$$
Certainly this bound is meaningful only if $n$ is large enough. This
provides
motivation for the study of curves over finite fields with many rational
points.

The purpose of this paper is to study {\it maximal curves}, that is,
curves $X$ over $\mathbb F_q$ whose number of rational points
$\#X(\mathbb F_q)$ reaches the Hasse-Weil upper bound. In this case one
knows that $q$ must be a square.

Let $k$ be the finite field with $q^2$ elements,
where $q$ is a power of a prime $p$. Let $X$ be a projective, connected,
non-singular
algebraic curve defined over $k$ which is maximal, that is, $\#X(k)$
satisfies
\begin{equation}\label{h}
\#X(k)=q^2 + 2gq + 1.
\end{equation}

Let $P\in X(k)$ and set $\d= g^{n+1}_{q+1}$ the $k$-linear system on
$X$ defined by the divisor $(q+1)P$. Then $n\ge 1$, and $\d$ is
independent of
$P$. In fact $\d$ is a simple base-point-free linear system on $X$
(Corollary \ref{R-Sti}, Remark \ref{rem-sem} (ii)). This allow us to apply
St\"ohr-Voloch's approach concerning Weierstrass point theory over finite
fields \cite{S-V}. Moreover, the dimension $n+1$ of $\d$ and the genus
$g$ are related by Castelnuovo's genus bound for curves in projective
spaces (\cite{C}, \cite[p.116]{ACGH}, \cite[Corollary 2.8]{Ra}).

It is known that $2g\le (q-1)q$
(\cite[V.3.3]{Sti}), and  that the Hermitian curve is the unique
maximal curve whose genus is $(q-1)q/2$
\cite{R-Sti}. Furthermore in \cite{F-T} we proved the following stronger
bounds for the genus, namely
$$
4g\le (q-1)^2\qquad {\rm or}\qquad 2g=(q-1)q.
$$
Moreover by using the already mentioned Castelnuovo's bound one can prove
that $4g>(q-1)^2$ if and only if $n=1$. Therefore, we assume from
now on that $n\ge 2$.

The Hermitian curve is a particular case of the following
type of curves. Let $m$ be a positive divisor of $q+1$, and let consider
\begin{equation*}
y^q + y = x^m.\tag{$\h$}
\end{equation*}
These curves are maximal (\cite[Thm. 1]{G-V}) and have very
remarkable properties (see e.g \cite{G-V}, \cite{Sch}).

Under a hypothesis on non-gaps at rational points we prove that
maximal curves are $k$-isomorphic to $\h$ for some $m\in \mathbb Z^+$.
\begin{theorem}\label{result1}
Let $X$ be a maximal curve of genus $g>0$. Assume
that there exists $P_0\in X(k)$ such that the first non-gap
$m_1$ at $P_0$ satisfies
$$
nm_1 \le q+1,
$$
where $n+1$ is the dimension of the complete linear system
defined by $(q+1)P_0$. Then one of the following possibilities is
satisfied
\begin{enumerate}
\item[(i)] $nm_1=q+1$, $2g=(m_1-1)(q-1)$, and $X$ is $k$-isomorphic to
${\mathcal H}_{m_1,q}$.
\item[(ii)] $nm_1=q$.
\end{enumerate}
\end{theorem}

{}From this theorem and a result due to Lewittes (see inequality \ref{le})
we obtain an analogous of
the main result in \cite{F-T}:
\begin{corollary}\label{result2}
Let $X$ be a curve satisfying the hypotheses of Theorem \ref{result1}.
Let $t\ge 1$ be an integer, and suppose that the genus $g$ of $X$ satisfies
$$
(q-1)(\frac{q+1}{t+1}-1)< 2g \le (q-1)(\frac{q+1}{t}-1).
$$

Then one of the following conditions is satisfied
\begin{enumerate}
\item[(i)] $t=n$, $2g=(q-1)(\frac{q+1}{t}-1)=(q-1)(m_1-1)$.
\item[(ii)] $t\ge n$, $2g\le
(q-1)(\frac{q+1}{n}-1)=(q-1)(\frac{m_1(q+1)}{q}-1)$.
\end{enumerate}
\end{corollary}
\begin{remarke}\label{conj}
In case $nm_1=q$ the authors actually conjecture that then $2g=(m_1-1)q$,
and $X$ is $k$-isomorphic to a curve whose plane model is given by
$F(y)=x^{q+1}$, where $F(y)$ is a $\mathbb F_p$-linear polynomial of
degree $m_1$. But we have not yet been able to prove
this. We notice that the veracity of this conjecture implies $t=n$ and
$2g=(\frac{q}{t}-1)q=(m_1-1)q$ in the statement (ii) of the above corollary.
\end{remarke}
\section{Preliminaries}
Throughout this paper we use the following notation:
\begin{itemize}
\item $k$ denotes the finite field with $q^2$ elements, where
$q$ is a power of a prime $p$. ${\bar k}$ denotes its algebraic closure.
\item By a curve we mean a projective, connected,
non-singular algebraic curve defined over $k$.
\item The symbol $X(k)$ (resp. $\kx$) stands for the set
of $k$-rational points (resp. the field of k-rational functions) of a curve
$X$.
\item If $x\in \kx$, ${\rm div}(x)$ (resp. $\div(x)$) denotes
the divisor (resp. the polar divisor) of $x$.
\item Let $P$ be a point of a curve. $v_P$ (resp. $H(P)$)) stands
for the valuation (resp. the Weierstrass semigroup) associated to $P$. We
denote by $m_i(P)$ the $i$th non-gap at $P$.
 \item Let $D$ be a divisor on $X$ and $P\in X$. We denote by
${\rm deg}(D)$ the degree of $D$, by $\supp(D)$ the support of $D$ and by
$v_P(D)$ the coefficient of $P$ in $D$. If $D$ is a $k$-divisor, we set
$$
L(D):= \{f\in \kx: {\rm div}(f)+D \succeq 0\},
$$
and $\ell(D):= {\rm dim}_k L(D)$. The symbol ``$\sim$" denotes
module linear equivalence.
\item The symbol $g^r_d$ stands for a linear system of dimension $r$
and degree $d$.
\end{itemize}
\subsection{Weierstrass points.}\label{wp} We summarize some results from
\cite{S-V}.
Let $X$ be a curve of genus $g$, $\d=g^r_d$ a base-point-free $k$-linear
system on $X$. Then associated to $P\in X$ we have the hermitian
$P$-invariants $j_0(P)=0<j_1(P)<\ldots<j_r(P)\le d$ of $\d$ (or simply
the $(\d,P)$-orders). This sequence is the same for all but finitely many
points. These finitely many points $P$, where exceptional $(\d,P)$-orders
occur, are called the $\d$-Weierstrass points of $X$. (If $\d$
is generated by a canonical divisor, we obtain the usual
Weierstrass points of $X$.) Associated to
$\d$ there exists a divisor $R$ supporting the $\d$-Weierstrass points. Let
$\epsilon_0<\epsilon_1<\ldots<\epsilon_r$ denote the $(\d,P)$-orders for
a generic $P\in X$. Then
\begin{equation}\label{ineq1}
\epsilon_i \le j_i(P),
\end{equation}
for $P\in X$, for each $i$, and
\begin{equation}\label{degR}
{\rm deg}(R)= (\epsilon_1+\ldots+\epsilon_r)(2g-2)+(r+1)d.
\end{equation}
Associated to $\d$ we also have a divisor $S$ whose support contains
$X(k)$. Its degree is
given by
\begin{equation*}
{\rm deg}(S)=(\nu_1+\ldots+\nu_{r-1})(2g-2)+(q^2+r)d,
\end{equation*}
where the $\nu_i's$ is a subsequence of the $\epsilon_i's$. More
precisely there exists an integer $I$ with $0<I\le r$ such that
$\nu_i=\epsilon_i$ for $i<I$ and $\nu_i=\epsilon_{i+1}$ otherwise.
Moreover for $P\in X(k)$ we have
\begin{equation}\label{ineq2}
v_P(S)\ge \sum_{i=1}^{r}(j_i(P)-\nu_{i-1}),
\end{equation}
and
\begin{equation}\label{ineq3}
\nu_i\le j_{i+1}(P)-j_1(P),
\end{equation}
for each $i$.
\subsection{Maximal curves}
Let $X$ be a maximal curve of genus $g$. In this
section we study  some arithmetical and  geometrical properties of $X$. To
begin with we have the following basic result which is containing in the
proof of \cite[Lemma 1]{R-Sti}. For the sake of completeness we state a
proof of it.
\begin{lemma1}\label{mult}
The Frobenius map $\pij$ (relative to $k$)
of the Jacobian $\j$ of $X$ acts
just as multiplication by $(-q)$ on $\j$.
\end{lemma1}
\begin{proof}
All the facts concerning Jacobians can be found in \cite[VI,
\S3]{L}. Let $\ell\not= p$ be a prime and let
$\T$ be the
Tate module of $\j$. Then the characteristic polynomial $P(\pij)(t)$ of
the
action $\pij$ on $\T$ is equal to $t^ {2g}L(1/t)$ where $L(t)$ denotes
the numerator of the Zeta function of $X$. Since $X$ satisfies
(\ref{h}), $L(t)=\prod_{i=1}^{2g}(1+qt)$ and thus $P(\pij)(t)= (t+q)^{2g}$.
Now we know that $\pij$ is diagonalizable \cite[Thm. 2]{Ta} and all its
eingenvalues are $-q$. This means that $\pij$ acts just as multiplication
by $-q$ on $\T$. Finally since the natural homomorphism of $\mathbb
Z$-algebras $$
{\rm End}(\j)\to {\rm End}(\T)
$$
is injective, the proof follows.
\end{proof}

Now fix $P_0\in X(k)$, and consider the map
$f=f^{P_0}: X\to \j$ given by $P\to [P-P_0]$. We have
$$
f\circ \pix = \pij\circ f,
$$
where $\pix$ denotes the Frobenius morphism of $X$ relative to $k$. Hence
from the above equality and Lemma \ref{mult} we get
\begin{corollary1}\label{frob}
$$
\pix(P)+qP \sim (q+1)P_0.
$$
\end{corollary1}

{}From this corollary it follows immediately the following:
\begin{corollary1}[[R-Sti, Lemma 1{]}]\label{R-Sti}
Let $P_0, P_1 \in X(k)$. Then $(q+1)P_1\sim (q+1)P_0$.
\end{corollary1}

Now, let consider the linear system $\d = g^{n+1}_{q+1}:=
|(q+1)P_0|$. Corollary \ref{R-Sti} says that $\d$ is a $k$-invariant of
the curve. In
particular its dimension $n+1$ is independent of $P\in X(k)$. Moreover
from Corollary \ref{R-Sti} we have that $q+1\in H(P_0)$ and hence
$\d$ is base-point-free. Consequently we can apply \cite{S-V} to $\d$.
\begin{thm1}\label{thm1} With notation as in \S\ref{wp} (for $\d$) we have:
\begin{enumerate}
\item[(i)] $\epsilon_{n+1}=\nu_n=q$;
\item[(ii)] $j_{n+1}(P)=q+1$ if $P\in X(k)$ and $j_{n+1}(P)=q$, otherwise;
\item[(iii)] $j_1(P)=1$ for each $P\in X$.
\end{enumerate}
\end{thm1}
\begin{proof} Statement (iii) for $P\in X(k)$ follows from (i), (ii) and
inequality (\ref{ineq3}). From Corollary \ref{frob} it follows (ii) and
$\epsilon_{n+1}=q$. Furthermore it also follows that $j_1(P)=1$ for
$P\not\in X(k)$: for let $P'\in X$ such that $\pix(P')=P$; then $P+qP'=
\pix(P')+ qP'\sim (q+1)P_0$.

Now we are going to prove that $\nu_n=\epsilon_{n+1}$. Let $P\in
X\setminus \{P_0\}$.
Corollary \ref{frob} says that $\pi(\pix(P))$ belongs to the osculating
hyperplane at $P$, where $\pi$ stands for the morphism associated to
$\d$. $\pi$ can be defined by a base $\{f_0,f_1,\ldots,f_{n+1}\}$ of
$L((q+1)P_0)$, where $v_P(f_i)\ge 0$ for each $i$. Let $x$ be a
separating variable of $k(X)\mid k$. Then by
\cite[Prop. 1.4(c), Corollary 1.3]{S-V}  the rational function
$$
w:= {\rm det} \begin{pmatrix}
f_0\circ\pix  & \ldots  &  f_{n+1}\circ\pix\\
D^{\epsilon_0}_x f_0 & \ldots & D^{\epsilon_0}_x f_{n+1}\\
\vdots      & \vdots &  \vdots\\
D^{\epsilon_n}_x f_0 & \ldots & D^{\epsilon_n}_x f_{n+1}
\end{pmatrix}
$$
satisfies $w(P)=0$ for each generic point $P$. Let $I$ be the smallest
integer such that the row $(f_0\circ\pix,\ldots,f_{n+1}\circ\pix)$ is a
linear combination of the vectors $(D^{\epsilon_i}_x
f_0,\ldots,D^{\epsilon_i}_x f_{n+1})$ with $i=0,\ldots,I$. Then according
to \cite[Prop. 2.1]{S-V} we find
$$
\{\nu_0<\ldots<\nu_n\}=\{\epsilon_0<\ldots<\epsilon_{I-1}
<\epsilon_{I+1}<\ldots<\epsilon_{n+1}\}.
$$
This concludes the proof.
\end{proof}
\begin{remark}\label{rem-sem} Let $X$ be a maximal curve.
\begin{enumerate}
\item[(i)] We claim that
$\nu_1=\epsilon_1=1$ (that is, the number $I$
in the above proof is bigger than one). In fact, suppose that $\nu_1>1$.
Then, due to the fact that $j_1(P)=1$ for each $P$, we can apply the proof
of
\cite[Thm. 1]{H-V} to conclude that $\#X(k)=(q+1)(q^2-1)-(2g-2)$. Then
from (\ref{h}) we must have $2g=(q-1)q$ and hence by \cite{F-T}, $n=1$, a
contradiction.
\item[(ii)] Let $P\in X(k)$. Due to Corollary \ref{R-Sti}
the fact that $j_1(P)=1$ and $j_{n+1}(P)=q+1$ is equivalent to have
$q,q+1\in H(P)$. This was noticed for some $P_0\in X(K)$ in \cite[Prop.
1]{Sti-X}.
\item[(iii)] Let $P\in X\setminus X(k)$. If $P \in X(\mathbb F_{q^4})$,
then $q-1\in H(P)$. If $P\not\in X(\mathbb F_{q^4})$, then $q\in H(P)$.

In fact, set  $i:= \min \{j\in \mathbb Z^+: \pix^j(P)=P\}$. Now applying
$\left(\pix^{i-1}\right)_*$ (see \cite[IV, Ex. 2.6]{Har}) to the
equivalence in Corollary \ref{frob} we get
$$
\pix(P)+(q-1)P \sim q\pix^{i-1}(P).
$$
Now the remarks follows from the fact that $\pix^{i-1}(P)=\pix(P)$
if and only if $P\in X(\mathbb F_{q^4})$.
\end{enumerate}
\end{remark}

In particular the above remark (ii) implies  that $\d$ is simple. Thus the
genus $g$ of $X$ and the dimension $n+1$ of $\d$ are related by
Castelnuovo's genus bound for curves in projective spaces (\cite{C},
\cite[p. 116]{ACGH}, \cite[Corollary 2.8]{Ra}). Thus
\begin{equation}\label{cast}
2g\le c(n,q):=M(q-n+e)\le
\begin{cases}
(2q-n)^2/4n   &  \text{if $n$ is even}\\
((2q-n)^2-1)/4n & \text{if $n$ is odd},
\end{cases}
\end{equation}
where $M$ is the biggest integer $\le q/n$ and $e:= q-Mn$.

We can also bound $g$ by using non-gaps at $P_0\in X(k)$. In fact,
Lewittes
\cite[Thm. 1(b)]{Le} proved that
$$
\#X(k)\le q^2m_1(P) +1,
$$
and hence from (\ref{h})
we conclude that
\begin{equation}\label{le}
2g\le q(m_1(P)-1).
\end{equation}
\begin{prop1}\label{ra1}
The following statements are equivalent:
\begin{enumerate}
\item[(i)] $\pi: X\rightarrow \mathbb P^{n+1}$ is a closed embedding,
i.e. $X$ is $k$-isomorphic to $\pi(X)$.
\item[(ii)] $\forall P \in X(\mathbb F_{q^4}): \pi(P)\in \mathbb
P^{n+1}(k) \Leftrightarrow P \in X(k)$.
\item[(iii)] $\forall P \in X(\mathbb F_{q^4}): q \in H(P)$.
\end{enumerate}
\end{prop1}
\begin{proof}
Let $P\in X$. Since $j_1(P)=1$ (cf. Theorem \ref{thm1} (iii)) we know
already that $\pi(X)$ is non-singular at all the branches centered at
$P$. Thus $\pi$ is an embedding if and only if $\pi$ is injective.
\begin{claim*}
If $\#\pi^{-1}(\pi(P))\ge 2$, then $P\in X(\mathbb F_{q^4})\setminus
X(k)$ and $\pi(P)\in \mathbb P^{n+1}(k)$.
\end{claim*}
\begin{proof} {\it (Claim).}
{}From  Corollary \ref{frob} it follows that $\pi^{-1}(\pi(P)\subseteq
\{P,\pi(P)\}$. Analogically we have $\pi^{-1}(\pi(\pix(P)))\subseteq
\{\pix(P), \pix^2(P)\}$. Thus if $\#\pi^{-1}(\pi(P))\ge 2$, then $P$
cannot be rational and $\pix^2(P)=P$, i.e. $P\in X(\mathbb
F_{q^4})\setminus X(k)$. Furthermore we have
$\pi(P)=\pi(\pix(P))=\pix(\pi(P))$, i.e. $\pi(P)\in \mathbb P^{n+1}(k)$.
\end{proof}
{}From this claim the equivalence (i) $\Leftrightarrow$ (ii) follows
immediately. As to the implication (i) $\Rightarrow$ (iii)
we know that $\dim |\pix(P)+qP-P-\pix(P)|=\dim |\pix(P)+qP| - 2$
(Corollary \ref{frob} and \cite[Prop.3.1(b)]{Har}), i.e. $q\in H(P)$.
Finally we want to conclude that $\pi$ is an embedding from
(iii). According to the above claim it is sufficient to show that
$\pi^{-1}(\pi(P))=\{P\}$ for $P \in X(\mathbb F_{q^4})$. Let $P\in X(\mathbb
F_{q^4})$. Because of Corollary \ref{frob} we know that
$\pi^{-1}(\pi(P))\subseteq \{P, \pix(P)\}$. Since $q\in H(P)$, there is a
divisor $D\in |qP|$ with $P\notin \supp(D)$. In particular
$$
\pix(P)+D\sim \pix(P)+qP\sim (q+1)P_0.
$$
Thus $\pi^{-1}(\pi(\pix(P)))\subseteq \supp(\pix(P)+D)$. So either
$\pi(P)\neq \pi(\pix(P))$ or $P=\pix(P)$. In both cases we have
$\pi^{-1}(\pi(P))=\{P\}$. This means altogether that $\pi$ is injective
and so indeed a closed embedding.
\end{proof}
\begin{prop1}\label{ra2}
Suppose that $\pi: X\rightarrow \mathbb P^{n+1}$ is a closed embedding.
Assume furthermore that there exist $r, s \in H(P_0)$ such that all
non-gaps at $P_0$ less than or equal to $q+1$ are generated by $r$ and
$s$. Then $H(P_0)$ is generated by $r$ and $s$. In particular the genus
of $X$ is equal to $(r-1)(s-1)/2$.
\end{prop1}
\begin{proof}
Let $x, y\in k(X)$ with $\div(x)=sP_0$ and $\div(y)=rP_0$. Since $q, q+1
\in H(P_0)$, the numbers $r$ and $s$ are coprime. Let  $\pi_2:
X\rightarrow \mathbb P^2$, $P\mapsto (1:x(P):y(P))$. Then the curves $X$
and $\pi_2(X)$ are birational and $\pi_2(X)$
is a plane curve given by an equation of type
$$
x^r+\beta y^s+\sum_{is+jr<nm} \alpha_{ij}x^iy^j=0,
$$
where $\beta,\alpha_{ij}\in k$ and $\beta\neq 0$. We are going to prove that
$\pi_2(P)$ is a non-singular point
of $\pi_2(X)$ for $P\neq P_0$. From this follows by \cite[Ch. 7]{Ful} that
$g=1/2(r-1)(s-1)$. Then by Jenkins \cite{J} we have $H(P_0)=\langle
r,s\rangle$.

Let $1,f_1,\ldots,f_{n+1}$ be a basis of $L((q+1)P_0)$, where
$n+1:=\dim |(q+1)P_0|$. Then there exist polynomials $F_i(T_1,T_2)\in
k[T_1,T_2]$ for $i=1,\ldots,n+1$ such that
$$
f_i=F_i(x,y)\ \mbox{on}\ X\qquad \mbox{for}\qquad i=1,\ldots, n+1.
$$
Consider the maps $\pi |(X\setminus\{P_0\}): X\setminus\{P_0\}\rightarrow
\mathbb  A^{n+1}$ given by $P\mapsto (f_1(P),\ldots,f_{n+1}(P))$;
$\pi_2 |(X\setminus\{P_0\}): X\setminus\{P_0\}\rightarrow \mathbb A^2$,
$P\mapsto (x(P),y(P))$; and $\phi :  \mathbb A^2\rightarrow \mathbb
A^{n+1}$, given by $(p_1,p_2)\mapsto
(F_1(p_1,p_2),\ldots,F_{n+1}(p_1,p_2))$.
Then the following diagram is commuting
$$
\Atriangle[X\setminus\{P_0\}`{\Bbb A}^2`{\Bbb A}^r;\pi_2`\pi`\phi].
$$
Thus we have for a point $P$ of
$X\setminus\{P_0\}$ and the corresponding local rings assigned to
$\pi(P), \pi_2(P)$ the commutative diagram
$$
\Atriangle[O_{\pi(X),\pi(P)}`O_{\pi_2(X),\pi_2(P)}`O_{X,P};f`c`h],
$$
where $h$ is injective since $k(X)=k(x,y)$, and $c$ is an isomorphism by
assumption. Thus $\pi_2{X}$ is non-singular at $\pi_2{P}$.
\end{proof}
\section{Proofs of Theorem 0.1 and Corollary 0.2}
Set $m:=m_1$. Recall that $n+1$ is by definition the dimension of
$\d:=|(q+1)P|$ for any $P\in X(k)$. Let $\pi$ be the morphism associated
to $\d$. By Remark \ref{rem-sem} (ii) we have $nm \ge q$, and hence by the
hypothesis on $m$ we get
$$
nm\in \{q,q+1\}.
$$
\subsection{Case: $nm=q+1$.}
\begin{prop1}\label{prop1} Let $X$ be a maximal curve of genus $g$.
Assume there exists $P_0\in X$ such that $nm_1(P_0)=q+1$. Then
$$
2g=(q-1)(m_1-1).
$$
\end{prop1}
\begin{proof} Since $m,q\in H(P_0)$ and $\gcd(m,q)=1$, then $2g\le
(m-1)(q-1)$ (see e.g. Jenkis \cite{J}). Now, $\pi$ can be
defined by $(1:y:\ldots:y^{n-1}:x:y^n)$ where $x, y \in k(X)$ such that
\begin{equation}\label{x,y}
\div(x)=qP_0\qquad {\rm and}\qquad \div(y)=mP_0.
\end{equation}
Let $P\in X\setminus \{P_0\}$. From the proof of \cite[Thm. 1.1]{S-V}, we
have that
\begin{equation}\label{vals}
v_P(y),\ldots,nv_P(y)
\end{equation}
are $(\d,P)$-orders. Thus by considering
a non-ramified point for $y:X\to \mathbb P^1$, and by
(\ref{ineq1}) we find
$$
\epsilon_i = i,\qquad {\rm for}\ \ i=1,\ldots, n.
$$
\begin{lemma1}\label{type-sem}
There are at most two types of $(\d,P)$-orders for $P\in X(k)$:
\begin{enumerate}
\item[(i):] $0,1,m,\ldots,(n-1)m,q+1$. Hence $w_1:=v_P(R)=
\frac{n((n-1)m-n-1)}{2}+2$.
\item[(ii):] $0,1,\ldots,n,q+1$. Hence $w_2:=v_P(R)=1$.
\end{enumerate}

Moreover, the set of the $\d$-Weierstrass points of $X$ coincides with
the set of $k$-rational points.
\end{lemma1}
\begin{proof} The statement on $v_P(R)$ follows from \cite[Thm. 1.5]{S-V}.
Let $P\in X(k)$. By Theorem \ref{thm1} we know that $1$ and $q+1$ are
$(\d,P)$-orders. We consider two cases:
\smallskip

(1) $v_P(y)=1$: With (\ref{vals}) this implies statement (ii).
\smallskip

(2) $v_P(y)>1$: From (\ref{vals}) it follows  $nv_P(y)=q+1$ and then we
obtain statement (i).
\smallskip

Let $P\in X\setminus X(k)$. By Theorem \ref{thm1} we have that
$j_{n+1}(P)=q$. If $v_P(y)>1$, then from (\ref{vals}) we get
$nv_P(x)=q=mn-1$
and hence $n=1$. Since by hypothesis $n>1$ then $v_P(y)=1$. This finish
the proof of the lemma.
\end{proof}

Let $T_1$ (resp. $T_2$)  denote the number of points $P\in X(k)$ whose
$(\d,P)$-orders are of type (i) (resp. type (ii)) in Lemma
\ref{type-sem}. Thus by (\ref{degR}) we have
$$
{\rm deg}(R)= (n(n+1)/2+q)(2g-2)+(n+2)(q+1)= w_1 T_1 + T_2,
$$
and by Riemann-Hurwitz applied to $y:X\to \mathbb P^1$
$$
2g-2=-2m + (m-1)T_1.
$$
Consequently, since $T_1+T_2 = \#X(k)= q^2+2gq+1$, from the above two
equations we obtain Proposition \ref{prop1}.
\end{proof}

Now we are going to prove the uniqueness part of the result. To begin
with we  generalize \cite[Lemma 5]{R-Sti}.
\begin{lemma1}\label{galois}
Let $X$ be a curve satisfying the hypotheses of Proposition \ref{prop1}.
Take
$y$ as in (\ref{x,y}). Then $k(X)\mid k(y)$ is a Galois cyclic extension.
\end{lemma1}
\begin{proof} Consider $y:X\to \mathbb
P^1(\bar k)$ as a map of degree $m=m_1$. From the proof of Lemma
\ref{type-sem} we see that $y$ has $(q+1)$
ramified points. Moreover, all of them are rational and totally ramified.
\begin{claim*} Let $P\in k\cup \{\infty\}$ such that $\#
y^{-1}(P)=m$. Then $y^{-1}(P)\subseteq X(k)$.
\end{claim*}
\begin{proof} {\it (Claim).} Let $P_1, \ldots, P_r\in k\cup \{\infty\}$
which are not ramified for $y$. Then $r\le q^2-q$. Let $n_i=\#
y^{-1}(P_i)\le m$. Since $2g=(q-1)(m-1)$ by Proposition \ref{prop1}, then
we have
$$
(q^2-q)m=\#X(k)-q-1= \sum_{i=1}^{r} n_i,
$$
from where it follows that $r=q^2-q$ and $n_i=m$ for each $i$.
\end{proof}
Now it follows that $k(X)\mid k(y)$ is Galois as in the proof of
\cite[Lemma 5]{R-Sti}. It is cyclic because there exists rational points
that are totally ramified for $y$.
\end{proof}
\begin{prop1}\label{prop2} Let $X$ be a curve as in Proposition
\ref{prop1}. Then $X$ is $k$ isomorphic to ${\mathcal H}_{m_1,q}$.
\end{prop1}
\begin{proof} Let $y$ be as in (\ref{x,y}).
\begin{claim}\label{eq} $X$ has a model plane given by an equation of type
$$
f(y)=v^m,
$$
where $f\in k[T]$ with ${\rm deg}(f)=q$, $f(0)=0$, and $v \in L(qP_0)$.
\end{claim}
\begin{proof} {\it (Claim \ref{eq}.)} We know that $k(X)\mid k(y)$ is
cyclic (Lemma \ref{galois}). Let $\sigma$ be a generator of $k(X)\mid
k(y)$. Set $V:= L(qP_0)$, $U:= L((n-1)mP_0)$. Then $\sigma\mid V\in {\rm
Aut}(V)$ and $\sigma\mid U = {\rm id}\mid U$. Since $p\nmid m$ we then
have that $\sigma\mid V$ is diagonalizable with an eigenvalue $\lambda$ a
primitive $m$-root of unity in $k$. Let $v\in V\setminus U$ be the
corresponding eigenvector for $\lambda$. Now since ${\rm Norm}_{k(X)\mid
k(y)}(v)= -v^m$ and since $v\in L(qP_0)$ we conclude the existence of
$f\in k[T]$ such that $f(y)=v^m$ and ${\rm deg}(f)=q$. Finally from the
fact that $y$ has exactly $(q+1)$ rational points as totally ramified
points, it follows that $f$ splits into linear factors in $k[T]$. Hence
we can assume $f(0)=0$.
\end{proof}

Now from the claim in the proof of Lemma \ref{galois}, Claim \ref{eq} and
$nm=q+1$ it follows that
$f^n(\alpha)-f^{nq}(\alpha)=0$ for $\alpha\in k$, and hence we obtain
\begin{equation*}
f^n(T) \equiv f^{nq}(T) \pmod{T^{q^2}-T}\tag{$*$}.
\end{equation*}
Set $f(T)=\sum_{i=1}^{q}a_iT^i$, $f^n(T)=\sum_{i=1}^{nq}b_i T^i$.
\begin{claim}\label{eq1}
$a_1\neq 0$, $a_i=0$ for $2\le i\le q-1$.
\end{claim}
\begin{proof} {\it (Claim \ref{eq1}).} $a_1\neq 0$ follows from $(*)$ and
$f(0)=0$. Suppose that $\{2\le i\le q-1: a_i\neq 0 \}\neq \emptyset$. Set
$$
t:={\rm min}\{1\le i\le q-1: a_i\neq 0\}\qquad \mbox{and} \qquad j:={\rm
max}\{1\le i\le q-1: a_i\neq 0\}.
$$
Due to the facts: multiplication by $q$ gives an automorphism of
$\mathbb Z/(q^2-1)\mathbb Z$, and $n-1+qt<q^2$ we then get
$b_{n-1+qj}=b^q_{(n-1)q+j}=na^{n-1}_q a_j\neq 0$. Then $nq={\rm deg}(f)\ge
n-1+qj$ implies together with $2n\le q+1$ that
\begin{equation*}
j+n-1<q \tag{$\dagger$}.
\end{equation*}

Then from $(*)$ we have
$b_{t+n-1}=na_ta^{n-1}_1\neq 0$.
Then again by $(*)$ and by $(\dagger)$ it follows that
$b_{q(t+n-1)}=b^q_{t+n-1}\neq 0$ which
implies $nq={\rm deg}(f^n)\ge q(t+n-1)$. But this contradicts to $t\ge
2$. Thus $a_i=0$ for $2\le i\le q-1$ and we are done.
\end{proof}
Write $f(T)=aT^q+bT$, $a,b\in k^*$. By Claim \ref{eq} we have that
$$
f(k)\subseteq \{\beta^m : \beta\in k\}=\cup{i=q}^{n-1} \xi^{im}\mathbb F_q,
$$
where $\xi$ is a primitive element of $k$. Now since $f(k)$ is a one
dimensional $k$-space, it follows that there exists $i\in\{0,\ldots,n-1\}$
such that $f(k)=\xi^{im}\mathbb F_q$. Set $x_1:= \xi^{-i}x$, $y_1:=
\epsilon y$, with $\epsilon$ being the unique element of $k^*$ such that
$$
{\rm Trace}^k_{\mathbb F_q}(\epsilon \alpha)=\xi^{-im}f(\alpha) \forall
\alpha \in k.
$$
These functions fulfil $y^q_1+y_1=x^m_1$ and we finish the proof of
Proposition \ref{prop1}.
\end{proof}

Now the proof of Theorem \ref{result1} (i) follows from the last two
propositions.
\begin{proof} {\it (Corollary \ref{result2}).} By Theorem \ref{result1}
we have two possibilities:
\smallskip

(1) $nm_1=q+1$: Then $2g=(m_1-1)(q-1)$ and statement (i) follows from the
hypothesis on $g$.
\smallskip

(2) $nm_1=q$: Here from $(q-1)(\frac{q+1}{t+1}-1)<2g$, (\ref{le}) and
$n<q$ we found $t\ge n$. The remaining part of (ii) follows from $2g\le
(q-1)(\frac{q+1}{t}-1)$.
\end{proof}
\begin{corollary1}\label{result3} Let $X$ be a maximal curve with genus
$g=(q-1)^2/4$. Then one of the following possibilities is satisfied:
\begin{enumerate}
\item[(i)] $X$ is $k$-isomorphic to ${\mathcal H}_{\frac{q+1}{2},q}$, or
\item[(ii)] For every point $P\in X(k)$ the first three non-gaps at $P$
are $\{q-1,q,q+1\}$.
\end{enumerate}
\end{corollary1}
\begin{proof}
{}From the hypothesis on $g$ and from (\ref{cast}) (applied to
$g^{n+1}_{q+1}$) we conclude $n+1\le 3$. Then by \cite{F-T} we must have
$n+1=3$. Let $P\in X(k)$ and let $m,q,q+1$ be the first three non-gaps.
Then $2m\ge q+1$ due to (\ref{le}). Moreover, $g\le g'$, where $g'$ is
the genus of the semigroup $\langle m,q,q+1\rangle $. We bound from above
$g'$ according to Selmer \cite[\S3.II]{Sel}. From that reference it
follows that $g'$ will be larger if $\gcd(m,q)=1$. So let assume this.
Define $s, t$ by $q+1=sq-tm$, $1<s<m$, $t>0$. Write $m=us+r$, $0\le r<s$.
Then we have ([p.6 loc. cit.])
$$
2g'= (m-1)(q-1) - ut(m-s+r).
$$
Hence by the hypothesis on $g$ we then have
$$
2ut(m-s+r)\le (q-1)(2m-q-1),
$$
and now it is easy to see that $2m=(q+1)$ or $m=q-1$. The first case
for some point $P\in X(k)$, and Theorem \ref{result1} (i) imply the
result.
\end{proof}
\subsection{Case: $nm_1=q$.}
As in the proof of Lemma \ref{type-sem} here one also has $\epsilon_i=\nu_i$
for $i=0,1\ldots, n-1$; $\epsilon_n=n$. However we cannot apply
\cite[Thm. 1.5]{S-V} to compute $v_P(R)$ for $P\in X(k)$.

If one can show that $\pi: X\rightarrow \mathbb P^{n+1}$ is a closed
embedding, then from Proposition \ref{ra2} we would have $2g=q(m_1(P)-1)$
for $P\in X(k)$.
\begin{remark} The hypothesis on the first non-gap of Theorem
\ref{result1} is necessary. In fact, consider the curve from Serre's list
(see \cite[\S4]{Se}) over $\mathbb F_{25}$, $g=3$. Then it is maximal. Let
$m,5,6$
be the first non-gaps at $P\in X(\mathbb F_{25})$. If $m=3$ from Theorem
\ref{result1} (i) we would have $g=4$. Thus $m=4$.
\end{remark}

\end{document}